\documentclass[submission,Phys]{SciPost}

\hypersetup{
    colorlinks,
    linkcolor={red!50!black},
    citecolor={blue!50!black},
    urlcolor={blue!80!black}
}

\usepackage[bitstream-charter]{mathdesign}
\DeclareSymbolFont{usualmathcal}{OMS}{cmsy}{m}{n}
\DeclareSymbolFontAlphabet{\mathcal}{usualmathcal}

\fancypagestyle{SPstyle}{
\fancyhf{}
\lhead{\colorbox{scipostblue}{\bf \color{white} ~SciPost Physics }}
\rhead{{\bf \color{scipostdeepblue} ~Submission }}

\fancyfoot[C]{\textbf{\thepage}}
}

\usepackage{paperstyle}

\begin{document}

\pagestyle{SPstyle}

\begin{center}{\Large \textbf{\color{scipostdeepblue}{
Fast elementwise operations on tensor trains with\\alternating cross interpolation
}}}\end{center}

\begin{center}\textbf{
Marc K.\ Ritter\,\orcidlink{0000-0002-2960-5471}
}\end{center}

\begin{center}
Center for Computational Quantum Physics, Flatiron Institute,\\ 162~5th~Avenue, New~York, NY~10010, USA
\\[\baselineskip]
\href{mailto:mritter@flatironinstitute.org}{\small mritter@flatironinstitute.org}
\end{center}

\section*{\color{scipostdeepblue}{Abstract}}
\textbf{\boldmath{%
Tensor trains (TTs), also known as matrix product states (MPS), are compressed representations of high-dimensional data that can be efficiently manipulated to perform calculations on the data. In many applications, such as TT-based solvers for nonlinear partial differential equations, the most expensive step is an elementwise multiplication or similar elementwise operation on multiple TTs.
Known error-controlled algorithms for such operations scale as \(\order(\chi^4)\), where \(\chi\) is the TT rank.
If the rank of the output is smaller than \(\chi^2\), it is possible to formulate algorithms with better scaling.
In this work, we present the alternating cross interpolation (ACI) algorithm that performs such operations in \(\order(\chi^3)\), while maintaining error control.
We demonstrate these properties on benchmark problems, achieving a significant speedup for TT ranks that are commonly encountered in practical applications.
}}

\vspace{\baselineskip}

\noindent\textcolor{white!90!black}{%
\fbox{\parbox{0.975\linewidth}{%
\textcolor{white!40!black}{\begin{tabular}{lr}%
  \begin{minipage}{0.6\textwidth}%
    {\small Copyright attribution to authors. \newline
    This work is a submission to SciPost Physics. \newline
    License information to appear upon publication. \newline
    Publication information to appear upon publication.}
  \end{minipage} & \begin{minipage}{0.4\textwidth}
    {\small Received Date \newline Accepted Date \newline Published Date}%
  \end{minipage}
\end{tabular}}
}}
}


\vspace{10pt}
\noindent\rule{\textwidth}{1pt}
\tableofcontents
\noindent\rule{\textwidth}{1pt}
\vspace{10pt}

\section{Introduction}

High-dimensional problems and problems involving a large range of length scales
are common challenges in computational physics,
as na\"ive approaches to such problems quickly exhaust the available computational resources.
A paradigmatic example is finding the ground state of an interacting many-body quantum system, a high-dimensional problem, which is commonly done using a compressed tensor network representation of the solution \cite{banuls_tensor_2023,orus_practical_2014}.
The most well-known of these methods is the density-matrix renormalization group (DMRG), which relies on  tensor trains (TTs), also known as matrix product states (MPSs) \cite{white_density_1992,mcculloch_density-matrix_2007,schollwock_density-matrix_2011}.
TT approaches have been used for many other high-dimensional problems, such as
evaluation of Feynman diagrams \cite{nunez_fernandez_learning_2022},
representing orbitals in quantum chemistry \cite{jolly_tensorized_2023},
and options pricing in financial mathematics \cite{arenstein_full_2025,glau_low-rank_2020,kastoryano_highly_2022,sakurai_learning_2025}.
Similarly, TTs can be used to solve problems of many length scales in a compressed format, such as
evaluating Bethe--Salpeter equations \cite{shinaoka_multiscale_2023,rohshap_two-particle_2025},
integrating over Brillouin zones \cite{ritter_quantics_2024},
calculating strain in 2D super-moir\'e materials \cite{fumega_correlated_2024},
and simulating fluid dynamics \cite{gourianov_quantum-inspired_2022,kornev_numerical_2023,peddinti_quantum-inspired_2024,gourianov_tensor_2025,holscher_quantum-inspired_2025}.

In many of these applications, the limiting factor is the \(\order(\chi^4)\) runtime scaling of known algorithms for multiplication and contraction of TTs \cite{stoudenmire_minimally_2010,chen_exponential_2018,ma_approximate_2024,camano_successive_2026}, where \(\chi\) is the rank of the input TT.
Multiplication of TTs is typically performed by contracting the TT with Kronecker-\(\delta\) tensor cores \cite{shinaoka_multiscale_2023,rohshap_two-particle_2025,bou-comas_quantics_2025}.
In the general case, where \(\chi'=\chi^2\), it is not possible to decrease the scaling below \(\order(\chi^4)\), since the resulting TT has \(\order(\chi'^2) = \order(\chi^4)\) elements.
In practical applications, such as differential equation solvers, the case \(\chi' \in \order(\chi)\) is, however, very common. For example, the input TT might correspond to the solution at a particular time step, and the result is the solution at the next time step. In this case, \(\chi\) is related to the correlations between different variables and length scales, which are governed by the differential equation to be solved. Therefore, they do not tend to increase infinitely at exponential speed. For \(\chi' \in \order(\chi)\), it is possible to formulate algorithms that scale as \(\order(\chi^3)\).
Such algorithms are particularly relevant for TT-based solvers for nonlinear differential equations, such as the Gross--Pitaevskii equation or the Navier--Stokes equations, as the cost of evaluating them is usually entirely dominated by the cost of evaluating elementwise products in the nonlinear terms \cite{gourianov_quantum-inspired_2022,kornev_numerical_2023,peddinti_quantum-inspired_2024,gourianov_tensor_2025,holscher_quantum-inspired_2025,bou-comas_quantics_2025,chen_solving_2026}.

This work presents an algorithm, called alternating cross interpolation (ACI), which computes elementwise operations, such as the Hadamard product, in an error-controlled way. It is closely related to the 2-site tensor cross interpolation (TCI) algorithm \cite{oseledets_tt-cross_2010,oseledets_tensor-train_2011,savostyanov_fast_2011,savostyanov_quasioptimality_2014,dolgov_parallel_2020,nunez_fernandez_learning_2022,nunez_fernandez_learning_2025}, combined with ideas from the alternating minimal energy method (AMEn) \cite{dolgov_alternating_2014}.
It is also inspired by the recursive sketched interpolation (RSI) algorithm presented in Ref.~\cite{meng_recursive_2026}, which scales as our ACI algorithm, but does not offer the same error control. To the best of our knowledge, other algorithms for efficient evaluation of elementwise operations that have been proposed in prior work \cite{michailidis_element-wise_2025, sun_hatt_2025, cazeaux_linear-scaling_2026} do not achieve scaling below \(\order(\chi^4)\) without imposing further assumptions.

\section{Alternating Cross Interpolation Algorithm}
\subsection{Problem statement}
Given a function
\(
    f: \mathbb{C}^N \rightarrow \mathbb{C},
\)
and \(N\) tensor trains
\begin{equation}
    x^n_{\sigma_1, \ldots, \sigma_{\seLL}} :=
    \raisebox{-16pt}{\includegraphics{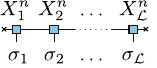}},
    \quad n = 1, \ldots, N,
\end{equation}
each with indices \((\sigma_1, \ldots, \sigma_{\seLL}) =: \bsigma \in \localset = \localset_1 \otimes \ldots \otimes \localset_\seLL\), we would like to obtain a tensor train \(y\) that approximates
\begin{equation}
    y_{\bsigma} :=
    \raisebox{-16pt}{\includegraphics{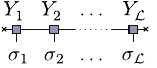}}
    \approx f(x^1_{\bsigma}, \ldots, x^N_{\bsigma})
    \label{eq:ydef}
\end{equation}
for all values of \(\bsigma\).
In practice, this means \(y\) should fulfill
\begin{equation}
    \left\|y - f(x^1, \ldots, x^N)\right\|_\infty =
    \max_{\bsigma\in\localset} \left|
        y_{\bsigma} - f(x^1_{\bsigma}, \ldots, x^N_{\bsigma})
    \right|
    \leq \tau
    \label{eq:ydef:maxnorm}
\end{equation}
for some user-specified tolerance \(\tau\).
Note that \(f\) operates \emph{element\-wise}: at each point \(\bsigma\), the output of \(f\) depends only on the value of \(x^1_\bsigma, \ldots, x^N_\bsigma\), all evaluated at the same point \(\bsigma\).
Common elementwise operations include addition of TTs
and the Hadamard product
\( x^1 \hadamard \cdots \hadamard x^N\).

\subsection{Overview}
ACI is an alternating optimization algorithm that sweeps across the tensor train, performing local updates at each pair of sites \(\ell, \ell+1\).
Efficient evaluation of the local updates is achieved by combining the CI-canonical gauge for the solution \(y\) with frame matrices for the inputs \(x\).
Each of these steps is explained below using tensor network diagrams \cite{orus_practical_2014}; a pseudocode in algebraic notation can be found in Appendix \ref{supp:pseudocode}.

\subsection{Alternating optimization}
Alternating optimization is a strategy to approach large optimization problems on TTs by casting the global optimization problem into a series of small local problems.
The idea is that successively solving each local problem eventually leads to global convergence.
The most famous example of such a method is the density matrix renormalization group (DMRG), which finds the ground state of a quantum system by solving a series of local eigenvalue problems \cite{schollwock_density-matrix_2011,mcculloch_density-matrix_2007,orus_practical_2014}.

The global problem considered in this work is to find tensors \(Y_1, \ldots, Y_\seLL\) that satisfy \Eq{eq:ydef:maxnorm}, ideally with quasi-optimal bond dimensions \(\chi'_1, \ldots, \chi'_N\). We approach this problem using 2-site alternating optimization, where we iterate over pairs of neighbouring sites. In each iteration, the algorithm solves the local problem of finding the optimal tensors \(Y_\ell\) and \(Y_{\ell+1}\), given current values of all other tensors \(Y_{\ell' \notin \{\ell, \ell+1\}}\). The solution \(y\) is then updated with the new \(Y_\ell\) and \(Y_{\ell+1}\) (\emph{local update}), before moving on to the next local problem.

The remainder of this section explains how to formulate a cheaply solvable local problem that leads to a solution of \Eq{eq:ydef:maxnorm}. To this end, we fix the gauge freedom in \(y\) to the CI-canonical gauge.

\subsection{CI-canonical gauge}
The CI-canonical gauge is defined through ordered index sets \(\Iset_\ell\) and \(\Jset_\ell\) \cite{nunez_fernandez_learning_2025}.
The left index sets \(\Iset_\ell \subseteq \localset_1 \otimes \ldots \otimes \localset_\ell\) contain multi-indices \(i = (i_1, \ldots, i_\ell)\) that can be used as indices to \(Y_1, \ldots, Y_\ell\).
Analogously, the right index sets \(\Jset_\ell \subseteq \localset_\ell \otimes \ldots \otimes \localset_\seLL\) contain multi-indices \(j = (j_\ell, \ldots, j_{\mathcal{L}})\) on \(Y_\ell, \ldots, Y_{\mathcal{L}}\).
Concatenating a multi-index \(i=(i_1, \ldots, i_\ell) \in \Iset_\ell\) with a multi-index \(j=(j_{\ell+1}, \ldots, j_\seLL) \in \Jset_{\ell+1}\) to \(ij=(i_1, \ldots, i_\ell, j_{\ell+1}, \ldots, j_\seLL)\) results in an index into the full tensor train \(y\).

We can now define \emph{slices} \(P_\ell\), \(T_\ell\) and \(\Pi_\ell\) of the tensor train \(y\) as
\begin{subequations}
\begin{align}
    \raisebox{-1.5pt}{\includegraphics{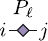}} \,&:= y_{ij},
    & i &\in \Iset_\ell, & j &\in \Jset_{\ell+1};\\
    \raisebox{-15pt}{\includegraphics{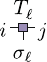}} \,&:= y_{i\sigma_\ell j},
    & i &\in \Iset_{\ell-1}, & j &\in \Jset_{\ell+1}; \\
    \raisebox{-15pt}{\includegraphics{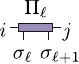}} \,&:= y_{i\sigma_{\ell}\sigma_{\ell+1} j},
    & i &\in \Iset_{\ell-1}, & j &\in \Jset_{\ell+2}\,.
\end{align}
    \label{eq:yslices}
\end{subequations}
In CI-canonical gauge, \(y\) is given by \cite{nunez_fernandez_learning_2025}
\begin{equation}
    y_{\sigma_1, \ldots, \sigma_N} =
    \raisebox{-16pt}{\includegraphics{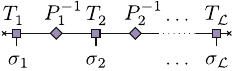}}
    \,.
    \label{eq:ycicanonical}
\end{equation}
The tensors \(Y_\ell\) in Eq.~\eqref{eq:ydef} can be recovered by multiplying all \(P_\ell^{-1}\) to the left or the right, e.g.\ \(Y_\ell = T_\ell P_\ell^{-1}\).

\subsection{Local problem}
We are now ready to formulate a local problem corresponding to \Eq{eq:ydef:maxnorm}.
We assume for the moment that good index sets \(\Iset_\ell, \Jset_\ell\) for \(y\) are known. Optimization of the index sets is discussed in a separate section below.
Given these index sets, each local tensor in \(y\) can be computed directly as
\begin{equation}
    \raisebox{-1.5pt}{\includegraphics{tndiagrams/Ptensor.pdf}} = y_{ij} = f(x^1_{ij}, \ldots, x^N_{ij})
    ,\quad\text{and}\quad
    \raisebox{-14pt}{\includegraphics{tndiagrams/Ttensor.pdf}} = y_{i\sigma_\ell j} = f(x^1_{i \sigma_\ell j}, \ldots, x^N_{i \sigma_\ell j}),
\label{eq:xslicesfory}
\end{equation}
combining Eq.~\eqref{eq:ydef} and Eq.~\eqref{eq:yslices}.
Thus, the local problem is reduced to computing the slices \(x^n_{ij}\) and \(x^n_{i \sigma_\ell j}\) on the index sets required for \(y\).
Evaluating a TT at a particular choice of indices costs \(\order(\chi^3 \eLL)\). Evaluating Eqs.~\eqref{eq:xslicesfory} na\"ively at each combination of \(i\), \(\sigma_\ell\), and \(j\) would therefore scale as \(\order(d\chi^3 \chi' \eLL N)\). This cost can be reduced by precomputing left and right \emph{frame matrices} \cite{dolgov_alternating_2014},
\begin{equation}
    \raisebox{-15pt}{\includegraphics{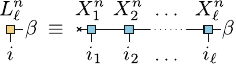}},
    \quad
    \raisebox{-15pt}{\includegraphics{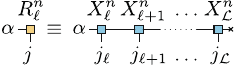}},
\end{equation}
with \(i = (i_1, \ldots, i_\ell) \in \Iset_\ell\) and \(j = (j_\ell, \ldots, j_{\eLL}) \in \Jset_\ell\).
There is an efficient method to compute \(L_\ell^n\) and \(R_\ell^n\) provided the index sets \(\Iset_\ell\) and \(\Jset_\ell\) satisfy the \emph{nesting conditions} \cite{nunez_fernandez_learning_2022,nunez_fernandez_learning_2025}:
The set \(\Iset_{\ell+1}\) is nested with respect to \(\Iset_{\ell}\) (denoted \(\Iset_{\ell} < \Iset_{\ell+1}\)) if \(\forall (i_1, \ldots, i_\ell, i_{\ell+1}) \in \Iset_{\ell+1}\),
removing the last index, \(i_{\ell+1}\), results in an element of \(\Iset_\ell\): \((i_1, \ldots, i_\ell) \in \Iset_\ell\).
Analogously, \(\Jset_{\ell-1} > \Jset_\ell\) if
\(
    \forall(j_{\ell-1}, j_\ell, \ldots, j_{\seLL}) \in \Jset_{\ell-1} :
    (j_{\ell}, \ldots, j_{\seLL}) \in \Jset_\ell
\).
If \(\Iset_\ell < \Iset_{\ell+1}\), \(L_{\ell+1}^n\) can be obtained from \(L_{\ell}^n\) by contracting
\begin{subequations}
\begin{gather}
    \raisebox{-15pt}{\includegraphics{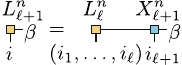}}
\intertext{
and truncating the index \(i\) such that \(i \in \Iset_{\ell+1}\). Analogously, \(R^n_{\ell-1}\) is obtained efficiently by contracting
}
    \raisebox{-15pt}{\includegraphics{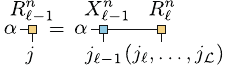}}
\end{gather}%
\label{eq:framesrecursive}%
\end{subequations}%
and truncating such that \(j \in \Jset_{\ell-1}\).
Given a set of frame matrices, the slices required for evaluating Eq.\ \eqref{eq:xslicesfory} are  evaluated in \(\order(\chi'\chi(\chi'+\chi)d)\) runtime by contracting
\begin{equation}
    x^n_{i \sigma_\ell j}
    =
    \raisebox{-15pt}{\includegraphics{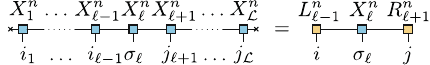}}.
\end{equation}
This way of efficiently evaluating the required slices of \(x^n\) using frame matrices is directly inspired by the AMEn method for solving linear systems using TT \cite{dolgov_alternating_2014}.
To our knowledge, the insight that frame matrices can be combined with the index sets of a CI-canonical form is novel. This idea is at the heart of the ACI algorithm.

\subsection{Optimizing index sets}
The remaining task of ACI is to find good index sets for \(y\).
We proceed exactly as in the 2-site TCI algorithm \cite{nunez_fernandez_learning_2025}.
We start with an initial guess \(\yinit\) for the solution; if no good initial guess is known, a random TT with bond dimensions equal to the minimum of the input bond dimensions is used.
\(\yinit\) is then put into CI-canonical form using the CI-canonicalization algorithm of Ref.\ \cite{nunez_fernandez_learning_2025}, which also gives nested initial right index sets \(\Jset_2, \ldots, \Jset_{\mathcal{L}}\). These index sets define initial right frame matrices \(R^n_2, \ldots, R^n_{\mathcal{L}}\). Left index sets \(\Iset_\ell\) and left frame matrices \(L^n_\ell\) are generated during the first sweep of ACI and do not need to be initialized, except for \(L^n_0 = [\,1\,]\).

We then sweep back and forth along the chain, starting with a left-to-right sweep. In each sweep, we perform local updates at each pair of sites \((\ell, \ell+1)\), which consist of the following steps. For each \(n \in \{1, \ldots, N\}\), the tensors \(\Pi^n_\ell\) are obtained by contracting
\begin{equation}
    \raisebox{-15pt}{\includegraphics{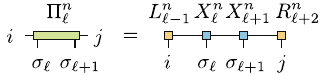}}.
    \label{eq:pifromframes}
\end{equation}
Now, the CI-canonical gauge (i) allows us to obtain \(\Pi_\ell\) efficiently by applying \(f\) elementwise to the tensors \(\Pi^n_\ell\). We then (ii) approximately factorize \(\Pi_\ell\) as
\begin{equation}
    \raisebox{-3pt}{\includegraphics{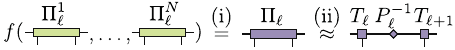}}.
    \label{eq:localupdate}
\end{equation}
The factorization into \(T_\ell P_\ell^{-1} T_{\ell+1}\) is obtained from a \emph{cross interpolation}, which chooses a subset of row and column indices of \(\Pi\) according to the \emph{maximum volume principle} \cite{goreinov_quasioptimality_2011,nunez_fernandez_learning_2025}. These row and column indices are then used to update \(\Iset_{\ell+1}\) and \(\Jset_{\ell}\), respectively.
During the sweep, these updates ensure that index sets are always nested with respect to site \(\ell\) \cite{nunez_fernandez_learning_2025}. The frame matrices on sites \(\ell, \ell+1\) are updated accordingly using Eq.~\eqref{eq:framesrecursive} after each local update.

\subsection{Complexity analysis}
A detailed analysis of runtime complexity can be found in Appendix \ref{supp:complexityanalysis}. Assuming \(\chi' = \chi\), i.e.\ identical bond dimensions of input and output TT, the most expensive step in each local update is evaluating Eq.~\eqref{eq:pifromframes}, which scales as \(\order(d^2\chi^3)\). Thus, the total runtime scales as \(\order( N_{\text{sweep}} \eLL N d^2 \chi^3)\).

\subsection{Alternatives}
The main alternative to ACI is the recently published recursive sketched interpolation (RSI) algorithm \cite{meng_recursive_2026}, whose runtime likewise scales as \(\order(\chi^3)\). The RSI algorithm is initialized by generating random sketch matrices of size \(\chi\times k\) for each site, where \(k\) is a user-specified parameter. It then performs a single sweep from left to right, contracting sketch matrices onto all tensors \(X^n_{\ell'>\ell}\) to the right of site \(\ell\).
This contraction results in four-leg tensors analogous to \(\Pi_\ell^n\), such that \(\Pi_\ell\) can be obtained by applying \(f\). \(\Pi_\ell\) is then factorized to obtain \(Y_\ell\). The main advantage of RSI over ACI is that the solution \(y\) is generated in a single left-to-right sweep. In exchange, the RSI algorithm cannot achieve the same error control that ACI provides through rank-adaptivity in subsequent sweeps.

In prior work, Hadamard products would usually be computed by attaching a Kronecker-\(\delta\) to each tensor in one of the inputs, and evaluating the product as an MPO-MPS contraction using one of the established algorithms \cite{stoudenmire_minimally_2010,chen_exponential_2018,ma_approximate_2024}. In the following examples, we compare the runtime of ACI to runtimes of a variational optimization algorithm for the MPO-MPS contraction implemented using the Julia packages TensorKit.jl and MPSKit.jl \cite{Devos_MPSKit_2026,devos_tensorkitjl_2025}.

\section{Examples}

\begin{figure}
    \centering
    \includegraphics{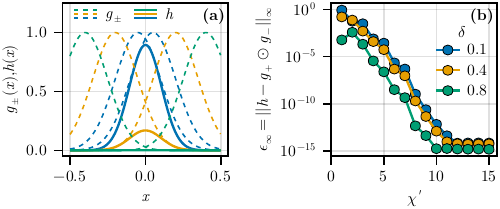}
    \caption{
        Multiplication of Gaussians for algorithm verification.
        \textbf{(a)}~The input functions \(g_\pm\) (dashed lines), given by Eq.~\eqref{eq:gaussians}, and their product \(h\), computed using ACI (solid lines).
        We set \(w = 0.15\) and vary \(\delta \in \{0.1, 0.4, 0.8\}\), discretize with \(\mathcal{L} = 25\) binary indices, and perform ACI with \(\chi' = 15\) on TTs generated with 2-site TCI.
        \textbf{(b)}~Maximum error between the exact product \(g_+ \odot g_-\) and the ACI output \(h\), as a function of output bond dimension \(\chi'\).
    }
    \label{fig:gaussians}
\end{figure}

\subsection{Gaussians}
As a simple example, we verify the ACI algorithm on the problem of multiplying two Gaussians,
\begin{equation}
    h(x) = g_+(x)\,g_-(x)
    \quad\text{with}\quad
    g_\pm(x) = \exp\left[- \frac{(x\pm\delta/2)^2}{2 w^2}\right].
    \label{eq:gaussians}
\end{equation}
The parameters \(\delta\) and \(w\) control the distance between peaks and the width of each peak, respectively. We discretize the interval \(x \in [x_{\min}, x_{\max}) = [-0.5, +0.5)\) using the quantics representation \cite{oseledets_approximation_2010}
\begin{equation}
   \frac{x - x_{\min}}{x_{\max} - x_{\min}}
   = \sum_{\ell=1}^{\mathcal{L}} 2^{-\ell} \sigma_\ell
   \label{eq:quanticsrep}
\end{equation}
with \(\mathcal{L} = 25\) binary indices \(\sigma_\ell\).
Quantics tensor trains (QTT) representing \(g_\pm\) are then obtained using 2-site TCI \cite{nunez_fernandez_learning_2025}, and used as input to ACI.

As shown in Fig.~\ref{fig:gaussians}(a), the peaks of \(f_\pm\) move apart with increasing \(\delta\), whereas \(h = g_+ \odot g_-\) always has a peak around \(x = 0\). For large \(\delta\), the optimal index sets for \(h\) are therefore very different from those for \(g_\pm\) \cite{meng_recursive_2026}.
This simple example thus verifies that ACI is capable of discovering structure in the solution which is not present in the inputs:
Fig.~\ref{fig:gaussians} shows that the maximum error between the ACI output and the exact solution,
\begin{equation}
    \epsilon_\infty =
    \left\| h - g_+ \odot g_- \right\|_\infty =
    \max_x \left| h(x) - g_+(x) g_-(x) \right|,
\end{equation}
decreases exponentially with increasing bond dimension \(\chi'\), down to numerical accuracy \(\epsilon_{\text{Float64}}\approx 10^{-14}\).
As the discretization grid contains too many points for explicit evaluation of \(\epsilon_\infty\), the error is approximated as the maximum over \(10^3\) random samples in this and the following examples.

\begin{figure}
    \centering
    \includegraphics{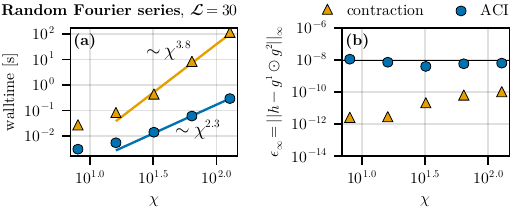}
    \caption{
    Hadamard product of random Fourier series, Eq.~\eqref{eq:fouriermultiplication}.
    Two functions of the form \eqref{eq:fourierseries} are represented as TT with \(\mathcal{L} = 30\) quantics indices, then multiplied elementwise with a tolerance of \(\tau = 10^{-8}\) using the ACI algorithm (this work, blue circles), and an algorithm based on MPO-MPS contraction (yellow triangles).
    \textbf{(a)}~Runtime needed to perform the multiplication.
    Lines are fitted power laws \(\sim\chi^p\), which are consistent with the theoretical scaling of \(\order(\chi^3)\) for ACI and \(\order(\chi^4)\) for the contraction-based Hadamard product.
    Runtimes were measured using a single thread on an AMD EPYC 9474F processor.
    \textbf{(b)}~Maximum error between the ACI output \(h\) and the exact product \(g^1 \odot g^2\). Both methods consistently reach the error tolerance \(\tau = 10^{-8}\) (black line).
    }
    \label{fig:fourier}
\end{figure}

\subsection{Random Fourier series}
To benchmark the runtime scaling and accuracy of the algorithm across many bond dimensions, we generate two random Fourier series
\begin{equation}
    g^{n}(x) = \sum_{k=0}^K \hat g^{n}_k e^{ikx}; \quad n = 1, 2
    \label{eq:fourierseries}
\end{equation}
with \(K+1\) complex Fourier coefficients \(\hat g_k^n\), whose real and imaginary parts are independently drawn from a uniform distribution on \([0, 1]\). The coefficients are then normalized such that \(\sum_k |\hat g^n_k|^2 = 1\).
We discretize these functions on the interval \(x \in [0, 1)\) using the quantics representation (Eq.~\eqref{eq:quanticsrep}) and obtain a QTT using 2-site TCI \cite{nunez_fernandez_learning_2025}.
Since the Fourier spectrum of \(f\) and \(g\) is band-limited with maximum frequency \(K\), the QTT have bond dimension \(\chi \in \order(\sqrt{K})\) \cite{lindsey_multiscale_2023}.
Multiplying them generates a new function,
\begin{equation}
    h(x) = g^1(x) \, g^2(x) =
    \sum_{k=0}^K \sum_{k'=0}^K \hat g^1_k \, \hat g^2_{k'} \, e^{i(k+k')x} =
    \sum_{q=0}^{2K} \left[ \sum_{q'=-q}^q \hat g^1_{(q+q')/2} \, \hat g^2_{(q-q')/2} \right] e^{iqx} =
    \sum_{q=0}^{2K} \hat h_q e^{iqx},
    \label{eq:fouriermultiplication}
\end{equation}
with up to \(2K+1\) coefficients \(\hat h_q\), which again has a bond dimension \(\chi' \in \order(\sqrt{K})\).

The runtime needed for the elementwise multiplication is shown in Fig.~\ref{fig:fourier}(a). Using ACI, it scales \(\sim\chi^{2.3}\), which is consistent with the theoretical asymptotic scaling of \(\order(\chi^3)\), accounting for the fact that \(\chi\) may not be large enough to reach the asymptotic regime.
For comparison, we compute the Hadamard product based on contracting with a Kronecker-\(\delta\) kernel.
We observe a scaling of runtime \(\sim \chi^{3.8}\), consistent with the theoretical scaling of \(\order(\chi^4)\).
Already at moderate bond dimensions of \(\chi \approx 100\), ACI achieves a speedup of a factor \( 10^2\) over the contraction algorithm.
As shown in Fig.~\ref{fig:fourier}(b), both methods reach the requested precision of \(\tau = 10^{-8}\) at all data points.

\begin{figure}
    \centering
    \includegraphics{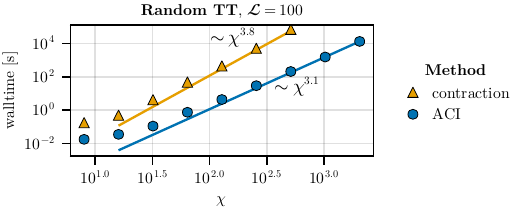}
    \caption{
    Similar to Fig.~\ref{fig:fourier}(a), but now for the Hadamard product of random TT. Instead of setting a fixed tolerance, the output bond dimension was set to \(\chi' = \chi\) here.
    }
    \label{fig:random}
\end{figure}

\subsection{Random TT}
To investigate the runtime scaling for larger \(\chi\), we multiply input TT that have been generated with random components drawn from a uniform distribution.
Except for normalization, this is equivalent to multiplying random wave functions.
The exact result of this multiplication has bond dimension \(\chi' = \chi^2\), but we apply the algorithms in a mode where we limit \(\chi'\) to \(\chi' \leq \chi\), instead opting to increase the truncation error. The resulting runtimes are shown in Fig.~\ref{fig:random}, and demonstrate scaling consistent with the theoretical \(\order(\chi^4)\) for the TT contraction algorithm and \(\order(\chi^3)\) for ACI.

\section{Conclusions and Outlook}
We have presented the alternating cross interpolation (ACI) algorithm for computing elementwise operations on tensor trains (TT). In benchmarks focused on computing Hadamard products of TT, we have demonstrated that the runtime scaling of ACI is \(\order(\chi^3)\) when the output bond dimension is \(\chi'\in\order(\chi)\), while maintaining error control through optimization of index sets.
In three examples, we have shown three key properties: ACI converges down to numerical accuracy provided the output bond dimension \(\chi'\) is large enough; ACI dynamically adjusts \(\chi'\) according to the specified error tolerance; and when \(\chi'\in\order(\chi)\), ACI scales as \(\order(\chi^3)\) over many orders of magnitude. The combination of these properties is unique to ACI.

These properties make ACI the method of choice for implementing nonlinearities in nonlinear differential equation solvers using TT.
In many nonlinear differential equations, all TT contractions are contractions with TT operators of constant, small bond dimension, such as the Fourier transform \cite{chen_quantum_2023,chen_direct_2026} or finite-difference derivative operators \cite{kazeev_low-rank_2012,gourianov_quantum-inspired_2022},
thus leaving elementwise operations necessary to evaluate nonlinearities as the only step scaling with \(\order(\chi^4)\).
For instance, TT-based solvers for computational fluid dynamics are limited by the evaluation of the nonlinear convection term \cite{gourianov_quantum-inspired_2022,gourianov_tensor_2025,holscher_quantum-inspired_2025}.
In these settings, switching to ACI with its improved scaling of \(\order(\chi^3)\) will result in improving the entire solver's scaling.
For further speedup, it is straightforward to use ACI within the patched elementwise multiplication scheme of Ref.~\cite{grosso_adaptive_2026}.

In other settings, such as evaluating equations arising in diagrammatic treatment of many-body theory, the presence of TT contractions scaling as \(\order(\chi^4)\) means that applying the ACI algorithm will improve the prefactor, but not the scaling of the total runtime \cite{rohshap_two-particle_2025}.
Nevertheless, improving prefactors is still desirable for practical applications, and using ACI removes the first of two bottlenecks in these applications, such that any improvement to TT contraction algorithms directly leads to speedup of TT-based solvers.

As for the algorithm itself, there are two obvious directions for future work.
Firstly, it is straightforward to generalize it to tree-shaped tensor networks in the same manner as TCI \cite{tindall_compressing_2024}.
Secondly, the structural similarity of ACI's local updates to RSI's sketching suggests that index sets used in CI-based algorithms can be understood as a particular type of sketch \cite{meng_recursive_2026}.
With sketching ideas being used in different settings to improve performance or reduce the runtime scaling of TT operations, a more complete understanding of the relation between different sketches and how CI fits into that framework may be fruitful.

\section*{Acknowledgments}
The author acknowledges useful discussions with J.\ von Delft, L.\ Devos, M.\ Frankenbach, O.\ Kovalska, M.\ Menon, I.\ V.\ Oseledets,  D.\ V.\ Savostyanov, and E.\ M.\ Stoudenmire.
The author thanks J.\ von Delft, L.\ Devos, M.\ Frankenbach, O.\ Kovalska, H.\ Shinaoka, and M.\ Wallerberger for feedback on the manuscript.
The author thanks L.\ Devos for providing a contraction-based Hadamard product implementation based on MPSKit.jl \cite{Devos_MPSKit_2026} and TensorKit.jl \cite{devos_tensorkitjl_2025}.

\paragraph{Code availability} All code and data necessary to reproduce the figures in this paper are available online \cite{coderepo}.
An implementation of the ACI algorithm is publicly available \cite{alternatingcrossinterpolation} as part of the tensor4all open-source libraries \cite{tensor4all.org}.

\paragraph{Funding information} The Flatiron Institute is a division of Simons Foundation.

\begin{appendix}
\numberwithin{equation}{section}

\section{Complexity analysis}\label{supp:complexityanalysis}
The three steps of each local update scale as follows.
\begin{enumerate}
    \item Contracting the frame matrices of size \(\chi' \times \chi\) with the input tensors \(X^n_\ell\) and \(X^n_{\ell+1}\) of size \(\chi\times d\times\chi\) costs \(\order(d\chi^2\chi'+d^2\chi{\chi'}^2)\) for each input \(x^n\) (Eq.\ \eqref{eq:pifromframes}).
    \item The function evaluation \(f(\Pi^1_\ell, \ldots, \Pi^N_\ell)\) costs \(\order(d^2{\chi'}^2)\) (Eq.\ \eqref{eq:localupdate}).
    \item Finally, factorizing \(\Pi_\ell\) to obtain index sets can be done in \(\order(d{\chi'}^3)\) (Eq.\ \eqref{eq:localupdate}).
\end{enumerate}
Each local update therefore scales as \(\order(Nd\chi^2\chi' + Nd^2\chi{\chi'}^2 + d{\chi'}^3)\). Since the local update is done \(\eLL N_{\text{sweep}}\) times, the total runtime cost is \(t \in\order(\eLL N_{\text{sweep}}[Nd\chi^2\chi' + d^2(N\chi+\chi'){\chi'}^2])\).
We can distinguish three common cases by the dependence of \(\chi'\) on \(\chi\):
\begin{enumerate}
    \item \(\chi' = \const \Rightarrow t \in\order(\eLL N_{\text{sweep}} N[d\chi^2 + d^2\chi])\).\\
    In this case, the output bond dimension is independent of \(\chi\).
    The contraction of the input tensors \(X^n_\ell\) to form frame matrices dominates the runtime (Eq.~\eqref{eq:framesrecursive}).
    \item \(\chi' \in \order(\chi) \Rightarrow t \in\order(\eLL N_{\text{sweep}} N d^2 \chi^3)\).\\
    The output and input are of comparable complexity, which is the case considered in the main text.
    The runtime is dominated by contracting frame matrices onto \(X^n_\ell\) to form the tensor \(\Pi_\ell^n\) in each local update (Eq.~\eqref{eq:pifromframes}).
    \item \(\chi' \in \order(\chi^2) \Rightarrow t \in \order(\eLL N_{\text{sweep}} d^2 [N\chi^5 + \chi^6])\).\\
    The output has maximum bond dimension. The runtime is dominated by the cost of factorizing the tensor \(\Pi_\ell\) (Eq.~\eqref{eq:localupdate}).
\end{enumerate}

In the common case where the elementwise operation to be computed is a Hadamard product, the most widely used alternative to ACI is a contraction with a diagonal tensor train, as described in e.g.\ Refs.\ \cite{shinaoka_multiscale_2023,rohshap_two-particle_2025}. The complexity of that algorithm is dominated by the cost of tensor train contraction (MPO-MPS contraction), for which algorithms with scaling \(\order(\eLL d^2\chi^4)\) are known \cite{stoudenmire_minimally_2010,chen_exponential_2018,ma_approximate_2024,camano_successive_2026}. ACI has better scaling compared to MPO-MPS contraction when \(\chi' \in \order(\chi^{4/3})\).

\section{Algorithms in pseudocode}
\label{supp:pseudocode}
The main alternating cross interpolation (ACI) algorithm is presented in Alg.~\ref{alg:elemmul}, with helper functions defined in Algs.~\ref{alg:ldu}, \ref{alg:ci}, \ref{alg:localupdate}, and \ref{alg:frame}.
For brevity, we use the following conventions in the pesudocode:
\begin{itemize}
    \item Index names: \(\sigma_\ell\) are always local (site) indices. \(i\) and \(j\) are reserved for multi-indices that are part of some index sets \(\Iset\) and \(\Jset\). All other indices without special properties are given Greek letters \(\alpha, \beta, \gamma\).
    \item Submatrices: To take submatrices and re-order rows and columns, we put index sets \(\Iset\) and \(\Jset\) in the indices of a matrix. For example, \(B = A_{\Iset, \Jset}\) means that \(B\) is the truncation of \(A\) to the rows in \(\Iset\) and columns in \(\Jset\).
    \item Index ranges: We use the notation \(m:n :=  (m, m+1, \ldots, n-1, n)\) to express index ranges.
    \item Indexed assignment: When indices appear on both sides of an assignment, this assignment is understood as assigning all elements. For example, \(A_{\alpha, 1}\gets B_{\alpha, 1}\) means
    \begin{algorithmic}
        \For{$\alpha \in \{1, \ldots, \operatorname{nrows}(A)\}$}
            \State \(A_{\alpha, 1} \gets B_{\alpha, 1}\)
        \EndFor
    \end{algorithmic}
    \item Einstein summation: Whenever an index appears on two or more elements that are being multiplied, that multiplication is to be understood as a contraction over that index. For example, \(C_{\alpha\beta} \gets A_{\alpha\gamma} B_{\gamma\beta}\) is equivalent to the matrix product \(C \gets AB\).
\end{itemize}
An implementation of ACI in the Julia programming language can be found online \cite{alternatingcrossinterpolation} as part of the tensor4all open-source libraries \cite{tensor4all.org}.

\begin{algorithm}
\caption{
The alternating cross interpolation (ACI) algorithm.
}
\label{alg:elemmul}
\begin{algorithmic}[1]
\Statex[-1]
{
\begin{tabular}[t]{lll}
\multicolumn{3}{l}{\functionname{AlternatingCrossInterpolation}}\\
\textbf{Input}~
& $f$ & elementwise function \\
& $X^1,\ldots, X^{N}$ & input tensor trains \\
& $Y^{\text{init}}$ & initial guess for $Y$; if nothing is known, initialize randomly.\\
& $\tau$ & absolute error tolerance \\
& $\hat\chi$ & maximum bond dimension \\
& $\variablename{maxiter}$ & number of iterations (sweeps) \\
\textbf{Output}~
& $Y$ & a tensor train approximating
$Y_{\bsigma} \approx f(X^1_{\bsigma}, \ldots, X^N_{\bsigma})$
\end{tabular}
}
\vspace{0.2\baselineskip}\hrule\vspace{0.2\baselineskip}
\For{$n \gets 1, \ldots, N$}
    \State{$L_0^{n} \gets \left[\,1\,\right]$}
    \Comment{Initialize first left frame matrix as \(1\times 1\) matrix.}
    \State $R^{n}_{\seLL+1} \gets \left[\,1\,\right]$
    \Comment{Analogous for last right frame matrix.}
\EndFor
\For{$\ell \gets \eLL, \eLL-1, \ldots, 2$}
    \Comment{Bring initial guess into CI-canonical form \cite{nunez_fernandez_learning_2025}. (Eq.~\eqref{eq:ycicanonical})}
    \State \(M_{\alpha, \sigma_\ell j} \gets [Y^{\text{init}}_\ell]_{\alpha, j}^{\sigma_\ell}\)
    \Comment{Reshape site tensor to matrix.}
    \State \(M', Y_\ell, \_, \Jset_\ell \gets \functionname{Crossinterpolate}(M, \tau, \hat\chi)\)
    \Statex\Comment{The row indices are not meaningful here, and are discarded.}
    \State \([Y_{\ell-1}]_{\alpha\beta}^\sigma \gets [Y_{\ell-1}]_{\alpha\gamma}^\sigma M'_{\gamma\beta}\)
    \Comment{Multiply left factor \(M'\) to the left.}
    \For{$n \gets 1, \ldots, N$}
    \Comment{Initialize right frames iteratively. (Alg.~\ref{alg:frame}, Eq.~\eqref{eq:framesrecursive})}
        \State $
            R^{n}_{\ell} \gets
            \functionname{Rightframe}(X^{n}_\ell, R^{n}_{\ell+1}, \indexset{J}_{\ell})
        $
    \EndFor
\EndFor
\For{$\variablename{iteration} \gets 1, \ldots, \variablename{maxiter}$}
\Comment{Main loop.}
\For{$\ell \gets 1, 2, \ldots, \eLL-1$}
    \Comment{Left-to-right sweep.}
    \State $
        [\Pi^{n}]^{\sigma_\ell, \sigma_{\ell+1}} \gets
        L^{n}_{\ell-1}
        [X^{n}_\ell]^{\sigma_\ell}
        [X^{n}_{\ell+1}]^{\sigma_{\ell+1}}
        R^{n}_{\ell+2}
    $
    \Statex\Comment{Assemble $\Pi^n$ in interpolative basis. (Eq.~\eqref{eq:pifromframes})}
    \State $
        \Pi^{\sigma_\ell, \sigma_{\ell+1}}_{\alpha, \beta} \gets
        f([\Pi^1]_{\alpha, \beta}^{\sigma_\ell, \sigma_{\ell+1}}, \ldots, [\Pi^N]_{\alpha, \beta}^{\sigma_\ell, \sigma_{\ell+1}})
    $
    \Statex\Comment{Apply elementwise operation $f$. (Eq.~\eqref{eq:localupdate})}
    \State{$
        Y_\ell, Y_{\ell+1}, \indexset{I}_\ell, \indexset{J}_{\ell+1}, \epsilon_\ell \gets \functionname{Localupdate}(\Pi, \tau, \hat\chi, \text{right})
    $}
    \Statex\Comment{Update bond $(\ell, \ell+1)$. (Alg.~\ref{alg:localupdate}, Eq.~\eqref{eq:localupdate})}
    \State{$
        L^n_\ell \gets \functionname{Leftframe}(X^n_\ell, L_{\ell-1}^n, \indexset{I}_\ell)
    $}
    \Comment{Update frame matrix. (Alg.~\ref{alg:frame}, Eq.~\eqref{eq:framesrecursive})}
\EndFor
\For{$\ell \gets \eLL-1, \ldots, 2, 1$}
    \Comment{Right-to-left sweep.}
    \State $
        [\Pi^{n}]^{\sigma_\ell, \sigma_{\ell+1}} \gets
        L^{n}_{\ell-1}
        [X^{n}_\ell]^{\sigma_\ell}
        [X^{n}_{\ell+1}]^{\sigma_{\ell+1}}
        R^{n}_{\ell+2}
    $
    \Statex\Comment{Assemble $\Pi^n$ in interpolative basis. (Eq.~\eqref{eq:pifromframes})}
    \State $
        \Pi^{\sigma_\ell, \sigma_{\ell+1}}_{\alpha, \beta} \gets
        f([\Pi^1]_{\alpha, \beta}^{\sigma_\ell, \sigma_{\ell+1}}, \ldots, [\Pi^N]_{\alpha, \beta}^{\sigma_\ell, \sigma_{\ell+1}})
    $
    \Statex\Comment{Apply elementwise operation $f$. (Eq.~\eqref{eq:localupdate})}
    \State $
        Y_\ell, Y_{\ell+1}, \indexset{I}_\ell, \indexset{J}_{\ell+1}, \epsilon_\ell \gets \functionname{Localupdate}(\Pi, \tau, \hat\chi, \text{left})
    $
    \Statex\Comment{Update bond $(\ell, \ell+1)$. (Alg.~\ref{alg:localupdate}, Eq.~\eqref{eq:localupdate})}
    \State $
        R^n_{\ell+1} \gets \functionname{Rightframe}(X^n_{\ell+1}, R_{\ell+2}^n, \indexset{J}_{\ell+1})
    $
    \Comment{Update frame matrix. (Alg.~\ref{alg:frame}, Eq.~\eqref{eq:framesrecursive})}
\EndFor
\If{$\chi_\ell$ have not increased $\And$ $\max_\ell \epsilon_\ell \leq \tau$}
    \State{\Break}
\EndIf
\EndFor
\end{algorithmic}
\end{algorithm}

\begin{algorithm}
\caption{In-place $LDU$ factorization of a matrix $M$ using the prrLU algorithm of Ref.~\cite{nunez_fernandez_learning_2025}.}
\label{alg:ldu}
\begin{algorithmic}[1]
\Statex[-1]
{
\begin{tabular}[t]{lll}
\multicolumn{3}{l}{\functionname{LDU}}\\
\textbf{Input}~
& $M$ & a $m \times n$ matrix to be factorized \\
& $\tau$ & truncation tolerance \\
& $\hat\chi$ & maximum bond dimension \\
\textbf{Output}~
& $L$ & a lower triangular matrix, \\
& $D$ & a diagonal matrix, and \\
& $U$ & an upper triangular matrix, such that $LDU \approx M$.\\
& $\Iset, \Jset$ & row and column index set \\
& $\epsilon$ & error estimate
\end{tabular}
}
\vspace{0.2\baselineskip}\hrule\vspace{0.2\baselineskip}
    \State $\Iset \gets (1, 2, \ldots, m)$
    \State $\Jset \gets (1, 2, \ldots, n)$
    \State $\chi \gets 1$
    \Repeat
        \State $i, j \gets \operatorname{argmax}_{\alpha\geq\chi, \beta\geq\chi} |M_{\alpha, \beta}|$
        \Statex\Comment{Find next pivot using a greedy version of the maxvol principle.}
        \State swap($M_{i, 1:n}, M_{\chi, 1:n}$) \Comment{Swap the next pivot to position \((\chi, \chi)\) in \(M\).}
        \State swap($M_{1:m, j}, M_{1:m, \chi}$)
        \State swap($\Iset_i, \Iset_\chi$) \Comment{Update \(\Iset\) and \(\Jset\) accordingly.}
        \State swap($\Jset_i, \Jset_\chi$)
        \State $M_{\chi+1:m, \chi+1:n} \!\gets\! M_{\chi+1:m,\chi+1:n} \!-\! M_{\chi:m, \chi} [M_{\chi, \chi}]^{-1} M_{\chi, \chi:n}$
        \Statex\Comment{Schur complement (Gaussian elimination step).}
        \State $\epsilon \gets |M_{\chi, \chi}|$ \Comment{Error estimate.}
    \Until{$\epsilon \leq \tau$ or $\chi = \hat\chi$}
    \State $L \gets \operatorname{lowertriangle}(M_{1:m,1:\chi})$
    \Comment{$L, D$, and $U$ are now in different parts of $M$.}
    \State $D \gets \operatorname{diagonal}(M_{1:\chi,1:\chi})$
    \State $U \gets \operatorname{uppertriangle}(M_{1:\chi,1:n})$
    \State $\Iset, \Jset \gets \Iset_{1:\chi}, \Jset_{1:\chi}$
\end{algorithmic}
\end{algorithm}

\begin{algorithm}
\caption{Cross interpolation (CI) of a matrix $M$ in a stable manner using $LDU$ factorization. Adapted from Ref.~\cite{nunez_fernandez_learning_2025}.}
\label{alg:ci}
\begin{algorithmic}[1]
\Statex[-1]
{
\begin{tabular}[t]{lll}
\multicolumn{3}{l}{\functionname{Crossinterpolate}}\\
\textbf{Input}~
& $M$ & a $m \times n$ matrix to be factorized \\
& $\tau$ & truncation tolerance \\
& $\hat\chi$ & maximum bond dimension \\
& \variablename{leftorright} & whether to multiply the factor \(P^{-1}\) to the left or the right factor \\
\textbf{Output}~
& $A$ & left factor, and \\
& $B$ & right factor, such that $AB \approx M$.\\
& $\Iset_{1:\chi}, \Jset_{1:\chi}$ & row and column index set \\
& $\epsilon$ & error estimate
\end{tabular}
}
\vspace{0.2\baselineskip}\hrule\vspace{0.2\baselineskip}
    \State $L, D, U, \indexset{I}, \indexset{J}, \epsilon \gets \functionname{LDU}(M, \tau, \hat\chi)$
    \Comment{(Alg.~\ref{alg:ldu})}
    \If{$\variablename{leftorright} = \text{left}$}
    \Comment{Left-orthogonal case.}
        \State $A \gets L \left[L_{\indexset{I},1:n}\right]^{-1}$
        \State $B \gets L_{\indexset{I},1:n} D U$
    \Else{} \Comment{Right-orthogonal case.}
        \State $A \gets L D U_{1:m, \indexset{J}}$
        \State $B \gets \left[U_{1:m, \indexset{J}} \right]^{-1} U$
    \EndIf
\end{algorithmic}
\end{algorithm}

\begin{algorithm}
\caption{Computes the factors \(Y_\ell, Y_{\ell+1}\) for the local update in AC (see Eq.~\eqref{eq:localupdate}).}
\label{alg:localupdate}
\begin{algorithmic}[1]
\Statex[-1]
{
\begin{tabular}[t]{lll}
\multicolumn{3}{l}{\functionname{Localupdate}}\\
\textbf{Input}~
& $\Pi$ & a two-site tensor (see Eq.~\eqref{eq:pifromframes}) \\
& $\tau$ & truncation tolerance \\
& $\hat\chi$ & maximum bond dimension \\
& \variablename{leftorright} & whether to multiply \(P^{-1}\) to the left or the right (see Eq.~\eqref{eq:yslices}) \\
\textbf{Output}~
& $Y, Y'$ & left and right factor, such that \(\Pi_{\alpha\beta}^{\sigma\sigma'} \approx Y_{\alpha\gamma}^{\sigma} Y'^{\sigma'}_{\gamma\beta}\) \\
& $\Iset_{1:\chi}, \Jset_{1:\chi}$ & row and column index set \\
& $\epsilon$ & error estimate
\end{tabular}
}
\vspace{0.2\baselineskip}\hrule\vspace{0.2\baselineskip}
    \State $M_{\alpha\otimes\sigma, \sigma' \otimes\beta} \gets \Pi_{\alpha\beta}^{\sigma, \sigma'}$
    \Comment{Reshape \(\Pi\) to a matrix.}
    \State $A, B, \indexset{I}, \indexset{J}, \epsilon \gets \functionname{Crossinterpolate}(M, \tau, \hat\chi, \variablename{leftorright})$
    \Comment{Factorize matricized $\Pi$. (Alg.~\ref{alg:ci})}
    \State $Y^{\sigma}_{\alpha\beta} \gets A_{\alpha\otimes\sigma, \beta}$
    \Comment{Reshape the matrix factors to tensor form.}
    \State $Y'^{\sigma'}_{\alpha\beta} \gets B_{\alpha, \sigma' \otimes\beta}$
\end{algorithmic}
\end{algorithm}

\begin{algorithm}
\caption{Iterative construction of left and right frame matrices (Eq.~\eqref{eq:framesrecursive}).}
\label{alg:frame}
\begin{minipage}[b]{0.48\linewidth}
\begin{algorithmic}[1]
\Statex[-1]
{
\begin{tabular}[t]{lll}
\multicolumn{3}{l}{\functionname{Leftframe}}\\
\textbf{Input}~
& $X$ & tensor for site \(\ell\) \\&& (size $\chi_L\times d \times \chi_R$) \\
& $L$ & left frame for sites \\&& \(1, \ldots, \ell-1\) \\
& $\Iset$ & index set for truncation \\&& of output frame \\
\textbf{Output}~
& $L'$ & left frame for sites \(1, \ldots, \ell\)
\end{tabular}
}
\vspace{0.2\baselineskip}\hrule\vspace{0.2\baselineskip}
    \State $L'_{i\sigma, \beta} \gets L_{i, \gamma} X^\sigma_{\gamma, \beta}$
    \Comment{\(i\sigma = (i_1, \ldots, i_{\ell-1}, \sigma)\)}
    \State $L' \gets L'_{\Iset, 1:\chi_R}$ \Comment{Truncate rows with $\indexset{I}$.}
\end{algorithmic}
\end{minipage}
\hfill
\begin{minipage}[b]{0.48\linewidth}
\begin{algorithmic}[1]
\Statex[-1]
{
\begin{tabular}[t]{lll}
\multicolumn{3}{l}{\functionname{Rightframe}}\\
\textbf{Input}~
& $X$ & tensor for site \(\ell\) \\&& (size $\chi_L\times d \times \chi_R$) \\
& $R$ & right frame for sites \\&& \(\ell+1, \ldots, \mathcal{L}\) \\
& $\Jset$ & index set for truncation \\&& of output frame \\
\textbf{Output}~
& $L'$ & left frame for sites \(1, \ldots, \ell\)
\end{tabular}
}
\vspace{0.2\baselineskip}\hrule\vspace{0.2\baselineskip}
    \State $R'_{\alpha, \sigma j} \gets X^\sigma_{\alpha, \gamma} R_{\gamma, j}$
    \Comment{\(\sigma j = (\sigma, j_{\ell+1}, \ldots, j_{\mathcal{L}})\)}
    \State $R' \gets R'_{1:\chi_L, \Jset}$ \Comment{Truncate columns with $\indexset{J}$.}
\end{algorithmic}
\end{minipage}
\end{algorithm}

\end{appendix}
\clearpage
\bibliography{ttmultiplication}


\end{document}